\documentclass[11pt]{article}
\usepackage{latexsym}
\usepackage{amsfonts}

\title{The independence number of a subset of an abelian group \\[.4in]}

\author{B\'{e}la Bajnok\thanks{Corresponding author} \\[.1in] Department of Mathematics, Gettysburg College \\
Gettysburg, PA 17325-1486 USA \\E-mail:  bbajnok@gettysburg.edu
\\[.2in] and
\\[.2in] Imre Ruzsa\thanks{Supported by Hungarian National Foundation for Scientific Research (OTKA), Grants No. T 25617, T 29759, T 38396} \\[.1in] Alfr\'ed R\'enyi Institute of Mathematics, Hungarian Academy of Sciences \\ Budapest, Pf. 127, H-1364 Hungary \\E-mail: ruzsa@renyi.hu }

\date{May 9, 2002}

\newtheorem{thm}{Theorem}[section]
\newtheorem{defin}[thm]{Definition}

\newtheorem{cor}[thm]{Corollary}
\newtheorem{prop}[thm]{Proposition}
\newtheorem{conj}[thm]{Conjecture}

\newcommand{\bc}[2]{{{#1}\choose{#2}}}

\begin{document}

\maketitle

\thispagestyle{empty}
\pagestyle{myheadings}
\markright{Independence number}

\newpage

\begin{abstract}

We call a subset $A$ of the (additive) abelian group $G$ {\it
$t$-independent} if for all non-negative integers $h$ and $k$ with $h+k \leq t$,
the sum of $h$ (not necessarily distinct) elements of $A$ does not equal
the sum of $k$ (not necessarily distinct) elements of $A$ unless $h=k$ and
the two sums contain the same terms in some order.  A {\it
weakly $t$-independent} set satisfies this property for sums of distinct terms.  We give some exact values and asymptotic bounds for the size of a largest $t$-independent set and weakly $t$-independent set
in abelian groups, particularly in the cyclic group ${\mathbb Z}_n$.

\newpage

\end{abstract}

\section{Introduction}

Our motivation for studying the independence number of subsets of an abelian group comes from spherical combinatorics.  It was shown by the first author in \cite{Baj:1998a} that if the set of integers $A=\{a_1,a_2,\dots, a_m\}$ forms a $3$-independent set in the cyclic group $\mathbb{Z}_n$ (as defined below), then the set of $n$ points $X=\{x_1,x_2,\dots, x_n\}$  with $$x_i=\frac{1}{\sqrt{m}} \cdot \left(\cos(\frac{2 \pi i a_1}{n}), \sin (\frac{2 \pi i a_1}{n}),  \dots, \cos (\frac{2 \pi i a_m}{n}), \sin (\frac{2 \pi i a_m}{n}) \right)$$ ( $i=1,2,\dots,n$) forms a {\em spherical $3$-design} on the unit sphere $S^{2m-1}$ (the case of $S^{2m}$ can be reduced to this case).  We believe that the concept of $t$-independence in $\mathbb{Z}_n$ and other abelian groups, extending some of the most well known concepts from additive number theory such as sum-free sets, Sidon sets, and $B_h$-sequences, is of independent interest; here we intend to provide the framework for a general discussion.

Throughout this paper $G$ denotes a finite abelian group with order $|G|=n \geq 2$, written in additive notation, and $A$ is a subset of $G$ of size $m \geq 1$.  For a positive integer $h$, we use the notation $$h \cdot A =\underbrace{A+A+ \cdots +A}_h = \{a_1+a_2 + \cdots + a_h | a_1,a_2,\dots,a_h \in A \}.$$  We introduce the following measure for the degree of independence of $A \subseteq G$.

\begin{defin} \label{t-independent}

Let $t$ be a non-negative integer and $A=\{a_1,a_2,\dots, a_m\}$.  We say that $A$ is a \emph{$t$-independent set} in $G$, if whenever $$\lambda_1a_1+\lambda_2a_2+\cdots +\lambda_m a_m=0$$ for some integers $\lambda_1, \lambda_2, \dots , \lambda_m$ with $$|\lambda_1|+|\lambda_2|+\cdots +|\lambda_m| \leq t,$$ we have $\lambda_1=\lambda_2= \cdots = \lambda_m=0$.  We call the largest $t$ for which $A$ is $t$-independent the \emph{independence number} of $A$ in $G$, and denote it by $\mathrm{ind}(A)$.  

\end{defin}

Equivalently, $A$ is a $t$-independent set in $G$, if for all non-negative integers $h$ and $k$ with $h+k \leq t$, the sum of $h$ (not necessarily distinct) elements of $A$ can only equal
the sum of $k$ (not necessarily distinct) elements of $A$ in a \emph{trivial} way, that is, $h=k$ and
the two sums contain the same terms in some order.  
We can break up our definition of $t$-independence into the following three requirements: 

\begin{equation} \label{one}
0 \not \in h \cdot A  \mbox{  for $1 \leq h \leq t$};
\end{equation}
\begin{equation} \label{two}
(h \cdot A) \cap (k \cdot A) = \emptyset  \mbox{  for $1 \leq h < k \leq t-h$};
\end{equation} and
\begin{equation} \label{three}
|h \cdot A | = \bc{m+h-1}{h}  \mbox{  for $1 \leq h \leq \left\lfloor \frac{t}{2} \right\rfloor$}.
\end{equation}
It is enough, in fact, to require these conditions for equations containing a total of $t$ or $t-1$ terms; therefore the total number of equations considered can be reduced to $2+(t-2)+1=t+1$.  

These conditions and their variations have been studied very vigorously for a long time; it might be worthwhile to briefly review some of the related classic and recent literature here.

Sets satisfying condition (\ref{one}) (or its versions where there is no limit on the number of terms and/or the terms need to be distinct -- see weak independence in Section 5) are called {\it zero-sum-free sets}.  For example, Erd\H{o}s and Heilbronn (\cite{ErdHei:1969a}, see also C15 in \cite{Guy:1994a}) asked for the largest number of distinct elements in the cyclic group $\mathbb{Z}_n$ so that no subset has sum zero.  Related results can be found in the papers of Alon and Dubiner \cite{AloDub:1993a}, Caro \cite{Car:1995a}, Gao and Hamidoune \cite{GaoHam:1998a}, Harcos and Ruzsa \cite{HarRuz:1997a}, and their references.

Sets satisfying the condition $(h \cdot A) \cap (k \cdot A) = \emptyset $ are called {\it $(h,k)$-sum-free sets}.  The first of these, $(1,2)$-sum-free sets, are simply called {\it sum-free sets}, and have a vast literature.  It is well known and not hard to see that, if $sf(G)$ denotes the maximum size of a sum-free set in $G$, then 
\begin{equation} \label{sumf} \frac{2}{7}n \leq sf(G) \leq \frac{1}{2}n, \end{equation}
where both inequalities are sharp as $sf(\mathbb{Z}_7)=2$ and $sf(\mathbb{Z}_2)=1$.
A comprehensive survey on sum-free sets in abelian groups can be found in Street's article \cite{WalStrWal:1972a}; see also the works of Erd\H{o}s \cite{Erd:1965a}, Alon and Kleitman \cite{AloKle:1990a}, and Cameron and Erd\H{o}s \cite{CamErd:1999a}.  Recently, $(h,k)$-sum-free sets have been investigated in cyclic groups of odd prime order by Bier and Chin \cite{BieChi:2001a}.  Other recent work on $(h,k)$-sum-free sets (among the positive integers rather than groups) includes a study of $(1,k)$-sum-free sets for $k \geq 3$ by Calkin and Taylor \cite{CalTay:1996a}, $(3,4)$-sum-free sets by Bilu \cite{Bil:1998a}, and $(h,k)$-sum-free sets by Schoen \cite{Sch:2000a}.  We return to sum-free sets in Section 2 as we study 3-independent sets.

Sets satisfying the condition that $h$-term sums of elements of $A$ be distinct up to the rearrangement of terms, as in (\ref{three}), are called {\it $B_h$-sequences}.  They have been studied very extensively among the positive integers, see the book of Halberstam and Roth \cite{HalRot:1983a}, sections C9 and C11 in Guy's book \cite{Guy:1994a}, and the survey paper of Graham \cite{Gra:1996a}.  The case $h=2$ is worth special mentioning; a $B_2$-sequence is called a {\it Sidon-sequence} after Sidon who introduced them to study Fourier series \cite{Sid:1932a}.  An excellent extensive survey of Sidon-sequences was written (in Hungarian) by Erd\H{o}s and Freud \cite{ErdFre:1991a}.  We review Sidon-sets and $B_h$-sequences as we study $t$-independence for $t \geq 4$ in Section 3.

There have recently been attempts to combine some of these conditions.  For example, Nathanson \cite{Nat:2000a} investigated sum-free Sidon sets in $\mathbb{Z}_n$; large Sidon sets were also used to construct small sum-free sets by Baltz, Schoen, and Srivastav \cite{BalSchSri:2000a}.  The present paper provides a more general setting for a discussion of these conditions.  

Zero-sum-free sets, $(h,k)$-sum-free sets, and $B_h$-sequences are only three of the many interesting families of linear equations among the integers or in an abelian group.  For a survey of general results and other known cases, see the papers of Ruzsa \cite{Ruz:1993a} and \cite{Ruz:1995a}, as well as sections C8 through C16, E12, E28, and E32 of Guy's wonderful book \cite{Guy:1994a}.  

We note that some, but not all, of the methods that we describe below can easily be modified to treat $t$-independence in non-abelian groups.  Sum-free (i.e. \emph{product-free}) sets and Sidon sets in non-abelian groups were discussed by Babai and S\'{o}s \cite{BabSos:1985a}; see also the recent results of Kedlaya in \cite{Ked:1997a} and \cite{Ked:1998a}.
 
Our objective in this paper is to study the size of a largest $t$-independent set in an abelian group $G$, which we denote by $s(G,t)$ (if $G$ has no $t$-independent subsets, we set $s(G,t)=0$).  

Since $0 \leq \mathrm{ind}(A) \leq n-1$ holds for every subset $A$ of $G$ (so no subset is ``completely'' independent), we see that $s(G,0)=n$ and $s(G,n)=0$.  It is also clear that $\mathrm{ind}(A)=0$ if and only if $0 \in A$, hence $s(G,1)=n-1.$  For the rest of the paper, we assume that $2 \leq t \leq n-1$. 

First we derive an upper bound for $s(G,t)$ using a simple counting argument, as follows.  Suppose that $A$ is a $t$-independent set in $G$ of size $m$.  Define $$\langle A, \lfloor t/2 \rfloor \rangle =\bigcup_{h=1}^{\lfloor t/2 \rfloor} h \cdot A.$$  Since $A$ is $t$-independent, by conditions (\ref{two}) and (\ref{three}) we see that $\langle A, \lfloor t/2 \rfloor \rangle$ has size exactly $$\sum_{h=1}^{\lfloor t/2 \rfloor} \bc{m+h-1}{h} = \bc{m+\lfloor t/2 \rfloor}{\lfloor t/2 \rfloor}-1.$$  

Therefore, the set $-\langle A, \lfloor t/2 \rfloor \rangle$, consisting of the negatives of the elements of $\langle A, \lfloor t/2 \rfloor \rangle$, also has this size; furthermore, we have $$\langle A, \lfloor t/2 \rfloor \rangle \cap -\langle A, \lfloor t/2 \rfloor \rangle = \emptyset.$$      
Additionally, by condition (\ref{one}), $$0 \not \in \langle A, \lfloor t/2 \rfloor \rangle \cup -\langle A, \lfloor t/2 \rfloor \rangle,$$ so $$n-1
\geq 2 \cdot |\langle A, \lfloor t/2 \rfloor \rangle|,$$  and therefore $$n \geq 2 \cdot \bc{m+\lfloor t/2 \rfloor}{\lfloor
t/2 \rfloor}-1.$$
Since for $t \geq 2$, $$2 \cdot \bc{m+\lfloor t/2 \rfloor}{\lfloor
t/2 \rfloor}-1 > 2 \cdot \frac{m^{\left\lfloor t/2 \right\rfloor}}{\left\lfloor t/2 \right\rfloor !},$$
we get the following result.

\begin{prop} \label{upperO} For every $t \geq 2$ we have $$s(G,t) <  \left( \frac{1}{2} \left\lfloor \frac{t}{2} \right\rfloor ! n \right)^{1/\left\lfloor t/2 \right\rfloor}.$$ \end{prop}

In section 3 we show that Proposition \ref{upperO} gives the correct magnitude of $s(G,t)$ if $G$ is the cyclic group $\mathbb{Z}_n$, that is $$s(\mathbb{Z}_n,t)=\Theta(n^{1/\lfloor t/2 \rfloor}).$$  Namely, we prove 

\begin{thm} \label{Zn-bound}

For every $\epsilon > 0$, $t \geq 2$, and large enough $n$, $$s(\mathbb{Z}_n,t) > \left(  \frac{1 - \epsilon}{t \cdot \left\lfloor (t+1)/2 \right\rfloor} \cdot n \right) ^{1/\left\lfloor t/2 \right\rfloor}.$$

\end{thm}

We have succeeded in finding the exact values of $s(\mathbb{Z}_n,t)$ only for $t \leq 3$.  We have already seen that $s(\mathbb{Z}_n, 1)=n-1$, and we can easily verify that \begin{equation} \label{Z_n2} s({\mathbb Z}_n,2)=\lfloor (n-1)/2 \rfloor \end{equation}   (the set $\{1,2,\dots, \lfloor (n-1)/2 \rfloor \}$ is 2-independent and has maximum size as we must have $A \cap -A = \emptyset$).  For $t=3$, in Section 2 we prove

\begin{thm} \label{3free}

$$s(\mathbb{Z}_n, 3) = \left\{
\begin{array}{cl}
\left\lfloor \frac{n}{4} \right\rfloor & \mbox{if $n$ is even}\\
\left(1+\frac{1}{p}\right) \frac{n}{6} & \mbox{if $n$ is odd, has prime divisors congruent to 5 $\pmod 6$,} \\ & \mbox{and $p$ is the smallest such divisor}\\
\left\lfloor \frac{n}{6} \right\rfloor & \mbox{otherwise}\\
\end{array}\right.$$

\end{thm}

Note that these results imply that the coefficient of $n$ in Theorem \ref{Zn-bound} cannot be improved for $t=2$ and $t=3$.  
 
Let us now turn to general abelian groups.  First note that, according to Proposition \ref{upperO} and Theorem \ref{Zn-bound}, the exponent of $n$ in the upper bound on $s(G,t)$ given in Proposition \ref{upperO} is sharp; for $$S(t)=\lim \sup \frac{s(G,t)^{{1/\left\lfloor t/2 \right\rfloor}}}{n}$$
we have $$ \frac{1 }{t \cdot \left\lfloor (t+1)/2 \right\rfloor}  \leq  S(t) \leq  \frac{1}{2} \left\lfloor \frac{t}{2} \right\rfloor ! .$$  These inequalities yield $S(2)=1/2$, and we later prove that $S(3)=1/4$.  We do not know the values of $S(t)$ for $t \geq 4$.

As for lower bounds, it is clear that we cannot expect a lower bound for $s(G,t)$ in terms of $n=|G|$ only; in fact, if the exponent $\kappa$ of $G$ is not more than $t$, then obviously $s(G,t)=0$.

We will use the following notations.  For a positive integer $h$, let the ``$h$-torsion'' subgroup of $G$ be \begin{equation} \label{Tor(G,h)}
\mathrm{Tor}(G,h)= \{x \in G | hx =0 \};
\end{equation}
then \begin{equation} \label{Ord(G,t)}
\mathrm{Ord}(G,h)= \cup_{h=1}^t  \mathrm{Tor}(G,h)
\end{equation}
is the set of those elements of $G$ which have order at most $t$.  By condition (\ref{one}), no element of $\mathrm{Ord}(G,h)$ can be in a $t$-independent set of $G$.    

We have already noted that $s(G,1)=n-1$, and we can now easily determine the value of $s(G,2)$: to get a maximum 2-independent set in $G$, take exactly one of each element or its negative in $G \setminus \mathrm{Ord}(G,2)$, hence we have

\begin{equation} \label{obs}  
s(G,2)=\frac{n-\mathrm{Ord}(G,2)}{2}.  
\end{equation}
As a special case, for the cyclic group of order $n$ we have (\ref{Z_n2}).

Note that if $\mathrm{Ord}(G,2)=G$ then $s(G,2)=0$; for $n \geq 2$ this occurs only for the elementary abelian 2-group.  If $\mathrm{Ord}(G,2) \not = G$ then, since $\mathrm{Ord}(G,2)$ is a subgroup of $G$, we have $1 \leq |\mathrm{Ord}(G,2)| \leq n/2$, and therefore we get the following.

\begin{prop}

If $G$ is isomorphic to the elementary abelian 2-group, then $s(G,2)=0$.  Otherwise $$\frac{1}{4}n \leq s(G,2) \leq \frac{1}{2}n.$$

\end{prop}

Let us now consider $t=3$.  As noted before, if $\mathrm{Ord}(G,3)=G$, then $s(G,3)=0$; this occurs if and only if $G$ is isomorphic to the elementary abelian $p$-group for $p=2$ or $p=3$.  In Section 2 we determine some exact values and sharp upper and lower bounds for $s(G,3)$ in terms of the exponent of $G$ (see Theorem \ref{3}); in particular, we prove the following.

\begin{thm} \label{3bounds}

If $G$ is isomorphic to the elementary abelian $p$-group for $p=2$ or $p=3$, then $s(G,3)=0$.  Otherwise $$\frac{1}{9}n \leq s(G,3) \leq \frac{1}{4}n.$$

\end{thm}
These bounds can be attained since $s(\mathbb{Z}_9,3)=1$ and $s(\mathbb{Z}_4,3)=1$.   

Studying $t$-independent sets in $G$ for larger values of $t$ seems considerably more difficult.  As a case in point, let $t \geq 4$, $\kappa >t$ be fixed, and consider the sequence of groups $$G_k=\mathbb{Z}_2^k \times \mathbb{Z}_{\kappa}$$ ($k=1,2,3 \dots$).  Suppose that $A$ is a $t$-independent set in $G_k$ and that for some $a_1,a_2 \in \mathbb{Z}_2^k$ and $x \in \mathbb{Z}_{\kappa}$, $g_1=(a_1,x)$ and $g_2=(a_2,x)$ are elements of $A$.  Then $g_1+g_1=g_2+g_2$ is a non-trivial equation in $G_k$, thus we conclude that $g_1=g_2$ and $|A| \leq \kappa$.  This implies that, as $k$ approaches infinity, we have $s(G_k,t)=O(1)$.

Therefore, in order to have a lower bound on $s(G,t)$ which tends to infinity with $n$, we must have more than just $\mathrm{Ord}(G,t) \not =G$ (as was the case for $t \leq 3$); although this will be sufficient for elementary abelian groups (see Corollary \ref{elem} below).  In general, we require that $|\mathrm{Ord}(G,t)|$ be not too large compared to $|G|=n$; in Section 4, we make this more precise by setting 
\begin{equation} \label{sigma}
\sigma(G,t)=\sum_{k=1}^t |\mathrm{Tor}(G,k)|
\end{equation}
and proving the following theorem.

\begin{thm} \label{general}

With $\sigma(G,t)$ as in (\ref{sigma}) above we have $$s(G,t) \geq \left\lfloor \left( \frac{n}{2 \sigma(G,t)} \right)^{1/t} \right\rfloor .$$

\end{thm}
 
Since obviously $\sigma(G,t) \leq t \cdot |\mathrm{Ord}(G,t)|$, we get the corollary that $$s(G,t) \geq \left\lfloor \left( \frac{n}{2 t \cdot |\mathrm{Ord}(G,t)|} \right)^{1/t} \right\rfloor.$$

It is worth comparing the lower bounds of Theorems \ref{Zn-bound} and \ref{general} for the cyclic group.  In $\mathbb{Z}_n$ we have
$$\sigma(\mathbb{Z}_n,t) = \sum_{h=1}^t |\mathrm{Tor}(\mathbb{Z}_n,h)| 
= \sum_{h=1}^t   \gcd(h,n)  \leq  \frac{t(t+1)}{2}.$$  Therefore both the exponent of $n$ and the coefficient are approximately twice as large in Theorem \ref{Zn-bound} as they are in Theorem \ref{general}.

Let us state a corollary of Theorem \ref{general} for elementary abelian $p$-groups.

\begin{cor} \label{elem}

Let $G$ be an elementary abelian $p$-group for a prime $p$.  If $p \leq t$, then $s(G,t)=0$; otherwise $$s(G,t) \geq \left(\frac{1}{2}n \right)^{1/t}.$$
\end{cor}

Finally, in Section 5, we examine {\em weak $t$-independence}: where for all non-negative integers $h$ and $k$ with $h+k \leq t$, sums of $h$ distinct elements of a set do not equal sums of $k$ distinct elements in a non-trivial way.  A comparison between sum-free and weak sum-free, as well as Sidon and weak Sidon sets, and $B_h$ sequences versus weak $B_h$ sequences, can be found in Ruzsa's papers \cite{Ruz:1993a} and \cite{Ruz:1995a}; there it was shown that their maximum sizes among the positive integers behave similarly.  This certainly does not hold for $t$-independence in abelian groups.  Denoting the maximum size of a weakly $t$-independent set in $G$ by $w(G,t)$, we find that for each fixed $t$, $\liminf w(G,t)$ tends to infinity with $n=|G|$ (see Theorem \ref{Weakly}), while obviously $\liminf s(G,t)=0$ for each $t \geq 2$.  Moreover, the weak independence number of each subset of $G$ can be arbitrarily large, even infinity, while its independence number cannot be more than $n$. 

This paper provides a modest attempt to discuss independence and weak independence in abelian groups.  Numerous interesting questions remain open, warranting further study.

\section{3-independent sets in abelian groups}

In this section we develop upper and lower bounds for $s(G,3)$, and provide exact values for some groups, including the cyclic group $\mathbb{Z}_n$.

Note that, by conditions (\ref{one}), (\ref{two}), and (\ref{three}), a subset $A$ of $G$ is 3-independent, if and only if $0 \not \in A, 0 \not \in A+A, 0 \not \in A+A+A$, and 
\begin{equation} \label{sumfr} (A+A) \cap A = \emptyset. \end{equation}  Sets satisfying (\ref{sumfr}) are called sum-free, and have been studied extensively; see \cite{WalStrWal:1972a} for a comprehensive survey.

It is not hard to determine bounds on the maximum size $sf(G)$ of a sum-free set in $G$.  First, note that the set $$\{\lfloor (n+1)/3 \rfloor , \lfloor (n+1)/3 \rfloor +1, \dots, 2 \lfloor (n+1)/3 \rfloor -1 \}$$ is a sum-free set in $\mathbb{Z}_n$; when $n$ is even, the larger set $$\{1,3,5, \dots, n-1 \}$$ is also sum-free.  Since $\lfloor (n+1)/3 \rfloor \geq \frac{2}{7}n$ when $n \geq 3$ is odd, we have $sf(\mathbb{Z}_n) \geq \frac{2}{7}n$.  Clearly, if $A$ is sum-free in $\mathbb{Z}_n$, then $H \times A$ is sum-free in $H \times \mathbb{Z}_n$ for any abelian group $H$, so the bound $sf(G) \geq \frac{2}{7}n$ holds for any $G$ (with $n>1$).  On the other hand, if $A$ is sum-free in $G$, then by (\ref{sumfr}) $n \geq |A+A|+|A| \geq 2|A|$, and we get the well known result that $$ \frac{2}{7}n \leq sf(G) \leq \frac{1}{2}n.$$
Note that these bounds are sharp, as $sf(\mathbb{Z}_7)=2$ and $sf(\mathbb{Z}_2)=1$.  

Our goal is to prove a similar result for 3-independent sets.  Let us start with the upper bound.

\begin{prop} \label{3upper}

Suppose that $A$ is a 3-independent set in $G$ of size $m$.
\begin{enumerate}

\item If none of the divisors of $n$ are congruent to $2 \pmod 3$, then $m \leq \frac{1}{6}n.$ 

\item Otherwise, let $p$ be the smallest divisor of $n$ which is congruent to 2 $\pmod 3$. Then we have
$m \leq \frac{1}{6}\left(1+\frac{1}{p}\right) n.$

\end{enumerate}
\end{prop}

{\it Proof.} Let $B= A \cup (-A)$. Since $A$ is 3-independent, we have $A \cap (-A) = \emptyset$, thus $ |B |=2m$. Furthermore, we have $$ B \cap (B+B) = \emptyset.$$
     
We apply Kneser's theorem \cite{Kne:1956a} to the set $B$.  It asserts that either we have $$ |B+B | \geq 2 |B |,$$ or there is a subgroup $H$ and an integer $k$ such that $B$ is contained
in $k$ cosets of $H$, and $B+B$ is equal to the union of $2k-1$ cosets.

In the first case we have $$ n \ge |B | + |B+B | \geq 3 |B | = 6m $$ and we are done.

Assume that the second possibility holds. Write $ |H |=d$ and $q=\frac{n}{d}$; we then have $$2m= |B | \leq dk  = \frac{n}{q} k,$$ hence 
\begin{equation} \label{m-bound} m \leq \frac{nk}{2q}. \end{equation}

Since $B$ cannot intersect any of the $2k-1$ cosets contained in $B+B$, we have $$ k + (2k-1) = 3k-1 \leq q;$$ if strict inequality holds here, then by (\ref{m-bound}) we have $$m \leq \frac{nk}{2(3k)}=\frac{1}{6}n,$$ and we are done again.

Otherwise, $q=3k-1  \equiv 2 \pmod{3}$, and from (\ref{m-bound}) we get $$ m \leq \frac{n}{2q} \frac{q+1}{3} = \frac{1}{6}\left(1+\frac{1}{q}\right) n \leq  \frac{1}{6}\left(1+\frac{1}{p}\right) n,$$ as claimed.  $\quad \Box$

Note that when $n$ is even, Proposition \ref{3upper} yields $s(G,3) \leq \frac{1}{4}n$.  This is clearly not sharp when $\mathrm{Ord}(G,2)$ is large; by (\ref{obs}) we must have $s(G,2)$ and therefore $s(G,3)$ small.  More precisely, we have the following upper bound.

\begin{prop} \label{3uppernew}

Let $\kappa$ be the exponent of $G$ and $\mathrm{Ord}(G,2)$ be the set of elements of $G$ with order at most 2.  If $\kappa$ is congruent to $2\pmod 4$, then $s(G,3) \leq \frac{1}{4}(n-|\mathrm{Ord}(G,2)|)$.

\end{prop} 

{\it Proof.}  When $\kappa \equiv 2 \pmod 4$, we can write $G=\mathbf{Z}_2^i \times G_1$, $ |G_1 |=n_1$ where $n_1$ is odd.  The case $n_1=1$ is obvious, so assume $n_1 \geq 3$.  Let $A$ be a 3-independent set in $G$ and write $|A|=m$.  We want to show that $m \leq 2^{i-2}(n_1-1)$.

Following the proof (and the notations) of Proposition \ref{3upper} above, we see that we either have $$m \leq \frac{1}{6}n,$$ or there is a subgroup $H$ of $G$ of index $q=3k-1$, such that $B=A \cup (-A)$ is contained in $k$ cosets of $H$, and $B+B$ is equal to the union of $2k-1$ cosets; in this case from (\ref{m-bound}) we have $$m \leq \frac{nk}{6k-2}.$$  Since $n_1 \geq 3$ implies $$\frac{1}{6}n \leq 2^{i-2}(n_1-1),$$ we can assume that the second possibility holds.  Furthermore, an easy computation shows that if $n_1 \geq 5$ and $k \geq 2$, then $$\frac{nk}{6k-2} \leq 2^{i-2}(n_1-1),$$ so we only need to consider
the cases of $n_1=3$ or $k=1$.

If $n_1=3$, then $G=\mathbf{Z}_2^i \times   \mathbf{Z}_3$. Let $B_0, B_1, B_2$ be the parts of $B$
in the three cosets of $\mathbf{Z}_2^i$. We have $B_0=\emptyset $, $B_2=-B_1$ and $ B_1 \cap  (B_2+B_2) = \emptyset $.  Now by an obvious pigeonhole argument we have $ |B_2 |\leq 2^{i-1}$, so $m \leq n/6$
again.  This concludes this case.

Finally assume $k=1$. Thus $H$ is a subgroup of index 2, $B$ is contained
in the coset $G \setminus H$ and $B+B=H$.  Since $\mathrm{Ord}(G,2)$ is a subgroup of $G$ of order $2^i$ and $H$ has order $2^{i-1}n_1$, we have $|H \cap \mathrm{Ord}(G,2)| \leq 2^{i-1}$, and therefore $|B| \leq n/2 - 2^{i-1}$, from which our claim follows.     $\quad \Box$

Now we turn to the lower bound.

\begin{prop} \label{3lower}

Let $\kappa$ be the exponent of $G$ and $\mathrm{Ord}(G,2)$ be the set of elements of $G$ with order at most 2.

\begin{enumerate}

\item  If $\kappa$ is divisible by 4, then $G$ contains a 3-independent set of size $\frac{1}{4}n$.

\item  If $\kappa$ is congruent to $2\pmod 4$, then $G$ contains a 3-independent set of size $\frac{1}{4}(n-|\mathrm{Ord}(G,2)|)$.

\item  Suppose that the odd positive integer $d$ divides $\kappa$.  Then $G$ contains a 3-independent set of size $\frac{1}{d} \lfloor \frac{d+1}{6} \rfloor n$.

\end{enumerate}

\end{prop}

{\it Proof.}  We construct the desired set explicitly as follows.  Suppose that the positive integer $d$ divides $\kappa$, and choose a subgroup $H$ and an element $g$ of $G$ so that $G$ is the union of the $d$ distinct cosets $H, H+g, \dots, H+(d-1)g$.  

Consider first the set $$A= \bigcup_{\frac{d}{6} < j < \frac{d}{3}} (H+jg).$$  It is clear that $A$ is 3-independent.  To determine the size of $A$, note that there are exactly $$f(d)=\left\lfloor \frac{d-1}{3} \right\rfloor - \left\lceil \frac{d+1}{6} \right\rceil +1$$ integers strictly between $\frac{d}{6}$ and $\frac{d}{3}$.  When $d$ is odd, $f(d)=\lfloor \frac{d+1}{6} \rfloor$, proving the last part of our Proposition.  Also, $f(4)=1$; so choosing $d=4$ when $\kappa$ is divisible by 4 proves the first case of our Proposition. 

Assume now that $\kappa$ is even, but not divisible by 4, and set $d=\kappa$.  Define $$A= \bigcup_{i=1}^{\lfloor \frac{\kappa}{4} \rfloor} (H+(2i-1)g).$$ It is easy to see that $A$ is 3-independent (note that $\kappa$ is even).  For the size of $A$ we have $$|A|=\left\lfloor \frac{\kappa}{4} \right\rfloor \cdot |H|=\frac{\kappa-2}{4} \cdot \frac{n}{\kappa}.$$  We add some elements to $A$ as follows.  Let $H'$ be a 2-independent set in $H$ of maximum size, and define $$A'=\{h'+\frac{\kappa}{2}g   |   h' \in H' \}.$$  Since $\frac{\kappa}{2}$ is odd, it is easy to check that $A \cup A'$ is 3-independent.  Using (\ref{obs}), we get $$|A'|=s(H,2)=\frac{\frac{n}{\kappa}-|\mathrm{Ord}(H,2)|}{2}.$$  But $G \cong H \times \mathbb{Z}_{\kappa}$, therefore $|\mathrm{Ord}(G,2)|=2 \cdot |\mathrm{Ord}(H,2)|$, and we get $$|A \cup A'|=\frac{\kappa-2}{4} \cdot \frac{n}{\kappa}+\frac{\frac{n}{\kappa}-|\mathrm{Ord}(H,2)|}{2}=\frac{1}{4}(n-|\mathrm{Ord}(G,2)|),$$ as claimed.     $\quad \Box$

We can now use Propositions \ref{3upper}, \ref{3uppernew}, and \ref{3lower} to establish the following bounds and exact values for $s(G,3)$.

\begin{thm} \label{3}

Let $\kappa$ be the exponent of $G$.  

\begin{enumerate}

\item

If $\kappa$ is divisible by 4, then $$s(G,3)=\frac{n}{4}.$$

\item

If $\kappa$ is even but not divisible by 4, then $$s(G,3)=\frac{1}{4}(n-|\mathrm{Ord}(G,2)|).$$

\item

If $\kappa$ (iff $n$) is odd and has prime divisors congruent to 5 $\pmod 6$, and $p$ is the smallest such divisor, then $$s(G,3)=\left(1+\frac{1}{p}\right) \frac{n}{6}.$$

\item

Finally, if $\kappa$ (iff $n$) is odd and has no prime divisors congruent to 5 $\pmod 6$, then $$\left\lfloor \frac{\kappa}{6} \right\rfloor \frac{n}{\kappa} \leq s(G,3) \leq \frac{n}{6}.$$

\end{enumerate}

\end{thm}

Since the exponent of the cyclic group $\mathbb{Z}_n$ is $n$, we have settled the value of $S(\mathbb{Z}_n, 3)$; see Theorem \ref{3free}.

In order to prove Theorem \ref{3bounds}, we should estimate $|\mathrm{Ord}(G,2)|$ when $\kappa$ is even.  Let $K$ be a cyclic subgroup of $G$ of order $\kappa$, and write $H=G/K$.  Then, since $\kappa$ is even, $|\mathrm{Ord}(G,2)|=2 \cdot |\mathrm{Ord}(H,2)|$.  Since obviously $|\mathrm{Ord}(H,2)| \leq |H|=\frac{n}{\kappa}$, we get $|\mathrm{Ord}(G,2)| \leq \frac{2n}{\kappa}$, and the first two cases of Theorem \ref{3} yield that when $\kappa$ is even, we have $$\left\lfloor \frac{\kappa}{4} \right\rfloor \frac{n}{\kappa} \leq s(G,3) \leq \frac{n}{4}.$$
So if $\kappa \geq 4$ is even, we get $$\frac{n}{6} \leq s(G,3) \leq \frac{n}{4},$$ and when $\kappa \geq 5$ is odd, Theorem \ref{3} implies $$\frac{n}{9} \leq s(G,3) \leq \frac{n}{5}.$$  In particular, we have proved Theorem \ref{3bounds}.

\section{$t$-independent sets in the cyclic group}

In this section we study the maximum size of a $t$-independent set in the cyclic group $\mathbb{Z}_n$.  In view of (\ref{Z_n2}) and Theorem \ref{3free}, it is enough to focus on $4 \leq t \leq n-1$.

For $t=4$ (in fact, for $t \geq 4$), condition (\ref{three}) requires that pairwise sums (or, equivalently, pairwise differences) of elements of $A$ be essentially distinct.  This condition has been studied extensively among the set of positive integers.  For a fixed positive integer $N$, a subset $B$ of $\{1,2, \dots, N\}$ with this property is called a {\it Sidon-sequence} after Sidon who introduced them to study Fourier series \cite{Sid:1932a}.  An excellent survey of Sidon-sequences was written (in Hungarian) by Erd\H{o}s and Freud \cite{ErdFre:1991a}.  Denoting the maximum cardinality of a Sidon-sequence in $\{1,2, \dots, N\}$ by $F_2(N)$, we have the classic result of Erd\H{o}s and Tur\'{a}n \cite{ErdTur:1941a} which says that for every $\epsilon >0$, $\delta >0$, and large enough $N$, 
\begin{equation} \label{F_2}
(1- \epsilon)\sqrt{N} < F_2(N) < (1+ \delta)\sqrt{N};
\end{equation} 
and it is a famous conjecture of Erd\H{o}s that, in fact, $|F_2(N)-\sqrt{N}|=O(1)$. 

More generally, a subset $B$ of $\{1,2, \dots, N\}$ with the property that all $h$-term sums of (not necessarily distinct) elements of $B$ are distinct, except for the order of the terms, is called a $B_h$-sequence.  Bose and Chowla \cite{BosCho:1962a} have shown that, if $F_h(N)$ denotes the maximum size of a $B_h$-sequence in $\{1,2, \dots, N\}$, then for every $\epsilon >0$ and large enough $N$, we have
\begin{equation} \label{F_h} F_h(N) > (1 - \epsilon)N^{\frac{1}{h}}.
\end{equation}
The simple counting argument leading to Proposition \ref{upperO}, noting that $h$-term sums of elements of $B$ are in the interval $[1,hN]$, yields the upper bound $$F_h(N) \leq (hh!N)^{\frac{1}{h}};$$ reducing the coefficient of $N$ is the subject of vigorous recent study (see C11 in \cite{Guy:1994a} and its references).  It is unknown whether the limit $$\lim_{n \rightarrow \infty} \frac{F_h(N)}{N^{1/h}}$$ exists for any $h \geq 3$.  For further references on $B_h$-sequences see the book by Halberstam and Roth \cite{HalRot:1983a}; sections C9 and C11 in Guy's book \cite{Guy:1994a}; and the survey paper of Graham \cite{Gra:1996a}.  

We use $B_h$-sequences to construct $t$-independent sets in the cyclic group $\mathbb{Z}_n$.  Elements of $\mathbb{Z}_n$ will be denoted by $0,1,2, \dots, n-1$.

\begin{prop} \label{bh}

Let $3 \leq t \leq n-1$ and \begin{equation} \label{N} N=\left\lfloor \frac{\left\lfloor n/t \right\rfloor}{\left\lfloor (t+1)/2 \right\rfloor} \right\rfloor, \end{equation} and suppose that $B$ is a $B_{\left\lfloor t/2 \right\rfloor}$-sequence in the interval $[1,N]$.  Then the set $$A=\left\{ \left\lfloor n/t \right\rfloor -b  |  b \in B \right\}$$ is $t$-independent in the cyclic group $\mathbb{Z}_n$.

\end{prop}

{\it Proof.}  First note that $3 \leq t \leq n-1$ guarantees that $N < \left\lfloor n/t \right\rfloor$, and thus, for each $a \in A$, we have $$0 <  \left\lfloor n/t \right\rfloor -N  \leq a \leq  \left\lfloor n/t \right\rfloor -1 < n/t.$$  We verify that requirements (\ref{one}), (\ref{two}), and (\ref{three}) hold.  

(i)  Let $1 \leq h \leq t$.  Since $g \in h \cdot A$ satisfies $$0 < h \cdot \left( \left\lfloor n/t \right\rfloor -N \right) \leq g \leq h \cdot \left(\left\lfloor n/t \right\rfloor -1 \right) <n,$$ we see that (\ref{one}) holds.

(ii)  To show (\ref{two}), let $1 \leq h < k \leq t-h$, and suppose, indirectly, that $g \in h \cdot A \cap k \cdot A$. Then we must have $$k \cdot \left( \left\lfloor n/t \right\rfloor -N \right) \leq g \leq h \cdot \left(\left\lfloor n/t \right\rfloor -1 \right).$$  Since $k \geq h+1$, this implies $$(h+1) \cdot \left( \left\lfloor n/t \right\rfloor -N \right) \leq h \cdot \left(\left\lfloor n/t \right\rfloor -1 \right),$$ or, equivalently, $$N \geq \frac{\left\lfloor n/t \right\rfloor +h}{h+1}.$$  But this contradicts (\ref{N}), since $ h \leq \left\lfloor (t-1)/2 \right\rfloor$.

(iii)  Finally, (\ref{three}) holds, since for $1 \leq h \leq \left\lfloor t/2 \right\rfloor$, $B$ is a $B_h$-sequence in $[1,N]$ and, using self-explanatory notation, $$a_1+a_2+ \cdots + a_h  = a_1'+a_2'+ \cdots + a_h' \mbox{  in  } \mathbb{Z}_n$$ implies $$a_1+a_2+ \cdots + a_h  = a_1'+a_2'+ \cdots + a_h' \mbox{  in  } \mathbb{Z}$$ and this further implies $$b_1+b_2+ \cdots + b_h  = b_1'+b_2'+ \cdots + b_h' \mbox{  in  } \mathbb{Z}.$$ $ \quad \Box$

Theorem \ref{Zn-bound} is now an easy corollary to Proposition \ref{bh} and (\ref{F_h}) (also using (\ref{Z_n2}) for $t=2$).

It is worthwhile to further analyze the cases $t=4$ and $t=5$ as follows.  Suppose that $A$ is a $t$-independent set in $\mathbb{Z}_n$ where $t \geq 4$.  Without loss of generality, we may assume that each $a \in A$ satisfies $1 \leq a \leq  \lfloor \frac{n-1}{2} \rfloor $ (replace $a$ by $n-a$, if necessary).  Note that $A$ is a Sidon-sequence in $\{1,2, \dots, \lfloor \frac{n-1}{2} \rfloor \}$.  Therefore, by (\ref{F_2}), we have $$|A| \leq (1+o(1)) \cdot \sqrt{n/2},$$ and this results in the following improvements.

\begin{cor} \label{s_n_lower_t=4,5}
For every $\epsilon >0$, $\delta >0$, and large enough $n$ we have $$\left(  \frac{1}{\sqrt{8}} - \epsilon  \right) \cdot \sqrt{n} \leq s(\mathbb{Z}_n,4) \leq \left(\frac{1}{\sqrt{2}}+\delta \right) \cdot \sqrt{n}.$$
$$\left(  \frac{1}{\sqrt{15}} - \epsilon  \right) \cdot \sqrt{n} \leq s(\mathbb{Z}_n,5) \leq \left(\frac{1}{\sqrt{2}}+ \delta \right) \cdot \sqrt{n}.$$

\end{cor}

We have determined the values of $s(\mathbb{Z}_n,4)$ and $s(\mathbb{Z}_n,5)$ for all $n \leq 200$.  While neither sequence is monotone, the values of $s(\mathbb{Z}_n,4)$ tend to grow more uniformly; we venture to state the following conjectures.

\begin{conj} \label{conj_t=4,5}  We have

\begin{enumerate}

\item

$\lim \frac{s(\mathbb{Z}_n,4)}{\sqrt{n}}= \frac{1}{\sqrt{3}},$ 

\item

$\lim \frac{s(\mathbb{Z}_n,5)}{\sqrt{n}} $ does not exist.

\end{enumerate}

\end{conj}

It is worth to recall, in comparison, that the sequence $s(\mathbb{Z}_n,2)/n$ is convergent while $s(\mathbb{Z}_n,3)/n$ is not.  A further observation on these sequences: by Theorem \ref{3free}, the sequence $s(\mathbb{Z}_n,3)$ is monotone for even values of $n$ (but not for the odd values); we find that, for $n \leq 200$, the sequence $s(\mathbb{Z}_n,5)$ is monotone for odd values of $n$ (but not for the even values).

\section{$t$-independent sets in abelian groups}

In this section we prove the general lower bound for $s(G,t)$ stated in Theorem \ref{general}.

Recall that for a positive integer $h$, we let the ``$h$-torsion'' subgroup of $G$ be $$
\mathrm{Tor}(G,h)= \{x \in G | hx =0 \};$$ and we also defined $$\sigma(G,t)=\sum_{h=1}^t |\mathrm{Tor}(G,h)|.$$

\begin{prop} \label{lower}

Suppose that $m$ is a positive integer for which $$n > \sigma(G,t) \cdot \bc{2m-2+t}{t}.$$  Then
$G$ has a $t$-independent set of size $m$.

\end{prop}

{\it Proof.} We use induction on $m$. For $m=1$ we have $n > \sigma(G,t)$, thus we can choose an element $$a \in G \setminus \bigcup_{h=1}^t \mathrm{Tor}(G,h).$$  Clearly, $\{a\}$ is then a $t$-independent set in $G$.

Assume now that our proposition holds for a positive integer $m$ and
suppose that $$n > \sigma(G,t) \cdot \bc{2m+t}{t}.$$ Since this value is greater than
$$ \sigma(G,t) \cdot \bc{2m-2+t}{t},$$ our inductive hypothesis implies that $G$ has a $t$-independent
set $A$ of size $m$.  

Define $B=A \cup (-A)$ and $$\langle B,t \rangle=\bigcup_{h=1}^t h \cdot B.$$

First, note that $\langle B,t \rangle$ has size at most $$\sum_{h=1}^t \bc{2m+h-1}{h} = \bc{2m+t}{t}-1.$$

For a fixed positive integer $h$ and group element $g$, define $$\mathrm{Root}_h(g)= \{x \in G | hx=g \}.$$  We prove the following

{\bf Claim.}  If $\mathrm{Root}_h(g) \not = \emptyset$, then $|\mathrm{Root}_h(g)|=|\mathrm{Tor}(G,h)|$.

{\it Proof of Claim.}  Fix $x \in \mathrm{Root}_h(g)$.  Our claim follows from the fact that $y \in \mathrm{Root}_h(g)$, if and only if, $y-x \in \mathrm{Tor}(G,h)$.

Now define $$C=\bigcup_{b \in \langle B,t \rangle} \bigcup_{h=1}^t \mathrm{Root}_h(b).$$  An obvious upper bound for the size of $C$ is $$ \sigma(G,t) \cdot \left( \bc{2m+t}{t}-1 \right).$$  Therefore, according to our inductive hypothesis, the set $G \setminus C$ is non-empty; fix $a \in G \setminus C$.

{\bf Claim.} $A \cup \{a\}$ is a $t$-independent set of size $m+1$ in $G$.

{\it Proof.}  To see that $A \cup \{a\}$ is of size $m+1$, note that $$a \in G \setminus C \subseteq G \setminus \langle B,t \rangle \subseteq G \setminus B \subseteq G \setminus A.$$

Now let $A=\{a_1,a_2,\dots,a_m\}$ and assume that $$\lambda_1a_1+\lambda_2a_2+\cdots +\lambda_m a_m+ \lambda a =0$$ for some integers $\lambda_1, \lambda_2, \dots , \lambda_m, \lambda$.  Suppose, indirectly, that $$1 \leq |\lambda_1|+|\lambda_2|+\cdots +|\lambda_m|+|\lambda| \leq t.$$  Set $$x= \lambda_1a_1+\lambda_2a_2+\cdots +\lambda_m a_m .$$ 

Note that $x \in \langle B,t \rangle$ and therefore $-x \in \langle B,t \rangle$; we also have $a \in \mathrm{Root}_{\lambda}(-x)$.

Without loss of generality, we can assume that $\lambda \geq 0$. If $1 \leq \lambda \leq t$, then we have $a \in \mathrm{Root}_{\lambda}(-x) \subseteq C$, a contradiction with the choice of $a$.  Otherwise, $\lambda =0$, from which $x=0$; a contradiction, since $A$ is $t$-independent.  This completes the proof of our claim and therefore our Proposition.  $ \quad \Box$

Now Theorem \ref{general} follows from Propositions \ref{lower}, by noting that 

\begin{eqnarray*}
\bc{2m-2+t}{t} & = & \prod_{k=1}^t \frac{2m-2+k}{k} \\
& = & \prod_{k=1}^t \frac{mk-(k-2)(m-1)}{k} \\
& \leq & (2m-1)m^{t-1} \\
& < & 2m^t.
\end{eqnarray*}

\section{Weakly $t$-independent sets in abelian groups}

When discussing solutions of equations in a set, it is natural to consider the version when we restrict ourselves to distinct solutions; this is referred to as the {\em weak} property.  A comparison between sum-free and weak sum-free, as well as Sidon and weak Sidon sets, and $B_h$ sequences versus weak $B_h$ sequences, can be found in Ruzsa's papers \cite{Ruz:1993a} and \cite{Ruz:1995a}; there it was shown that their maximum sizes among the positive integers behave similarly.  This certainly does not hold for $t$-independence in abelian groups, as we see in this section.

For a positive integer $h$, we use the notation $$h \star A =\{a_1+a_2 + \cdots + a_h | a_1,a_2,\dots,a_h \in A \mbox{  are    distinct}  \}.$$  We introduce the following measure for the degree of weak independence of $A \subseteq G$.

\begin{defin} \label{weak t-independent}

Let $t$ be a non-negative integer and $A=\{a_1,a_2,\dots, a_m\}$.  We say that $A$ is a \emph{weakly $t$-independent set} in $G$, if whenever $$\lambda_1a_1+\lambda_2a_2+\cdots +\lambda_m a_m=0$$ for some integers $\lambda_1, \lambda_2, \dots , \lambda_m \in \{-1,0,1\}$ with $$|\lambda_1|+|\lambda_2|+\cdots +|\lambda_m| \leq t,$$ we have $\lambda_1=\lambda_2= \cdots= \lambda_m=0$.  We call the largest $t$ for which $A$ is weakly $t$-independent the \emph{weak independence number} of $A$ in $G$ and denote it by $\mathrm{wind}(A)$; if $A$ is weakly $t$-independent for every $t$, we set $\mathrm{wind}(A)=\infty$.  

\end{defin}

Equivalently, $A$ is a weakly $t$-independent set in $G$, if for all non-negative integers $h$ and $k$ with $h+k \leq t$, the sum of $h$ distinct elements of $A$ can only equal
the sum of $k$ distinct elements of $A$ in a \emph{trivial} way, that is, $h=k$ and
the two sums contain the same terms in some order.  

This time, we have the following three requirements: 

\begin{equation} \label{onestar}
0 \not \in h \star A  \mbox{  for $1 \leq h \leq t$};
\end{equation}
\begin{equation} \label{twostar}
(h \star A) \cap (k \star A) = \emptyset  \mbox{  for $1 \leq h < k \leq t-h$};
\end{equation} and
\begin{equation} \label{threestar}
|h \star A | = \bc{m}{h}  \mbox{  for $1 \leq h \leq \left\lfloor \frac{t}{2} \right\rfloor$}.
\end{equation}

The difference between independence and weak independence can be illustrated by the following examples in $G=\mathbb{Z}_{30}$: we see that $\mathrm{ind}(\{1,2,4,8,16\})=2$ (as $1+1=2$), and $\mathrm{wind}(\{1,2,4,8,16\})=3$ (we have $2+4+8+16=0$); but $\mathrm{ind}(\{1,2,4,8\})=2$ still, yet $\mathrm{wind}(\{1,2,4,8\})=\infty$.

For a non-negative integer $t$, we let $w(G,t)$ denote the size of a maximum weakly $t$-independent set in $G$; we also set $w(G,\infty)$ to be the largest size of a subset $A$ of $G$ for which $\mathrm{wind}(A)=\infty$.  It is easy to see that $$w(G,0)=n,$$  $$w(G,1)=n-1,$$ and $$w(G,2)=\frac{n+|\mathrm{Ord}(2)|-2}{2};$$ the last equation results from the fact that no element of $G$, other than those of order 2, can be in a weakly 2-independent set together with its negative.  On the other end, if the invariant factor decomposition of $G$ contains $s$ terms, then $w(G,\infty)\geq s$; in particular, if $\kappa$ denotes the exponent of $G$, then $$w(G,t) \geq w(G,\infty) \geq \frac{\log n}{\log \kappa}$$ holds for every $t$.  

It is not hard to prove the following stronger result.

\begin{thm} \label{Weakly}  For $t \geq 2$ we have
$$ \left( \frac{t!}{2^t} n \right)^{1/t}- \frac{t}{2} < w(G,t) <  \left( \left\lfloor \frac{t}{2} \right\rfloor ! n \right)^{1/\left\lfloor t/2 \right\rfloor}+\frac{t}{2}.$$
\end{thm}

{\it Proof.}  Let us first prove two claims.

{\bf Claim 1.}  Suppose that $m$ is a positive integer for which $$n > \sum_{h=1}^t \bc{2m-2}{h}+1.$$  Then
$G$ has a weakly $t$-independent set of size $m$.

{\it Proof of Claim 1.} The proof will be similar (but simpler than) that of Proposition \ref{lower}.  We use induction on $m$. For $m=1$ we have $n \geq 2$, and 
clearly $\{a\}$ is a weakly $t$-independent set in $G$ whenever $a  \not = 0$.

Assume now that our proposition holds for a positive integer $m$ and
suppose that $$n > \sum_{h=1}^t \bc{2m}{h}+1.$$ Since this value is greater than
$$ \sum_{h=1}^t \bc{2m-2}{h}+1,$$ our inductive hypothesis implies that $G$ has a weakly $t$-independent
set $A$ of size $m$.  

Define $B=A \cup (-A)$ and $$\langle B, t \rangle ^*=\bigcup_{h=1}^t h \star B.$$

Then $|\langle B, t \rangle ^*| \leq \sum_{h=1}^t \bc{2m}{h}$.  Therefore, we can choose an $a \in G \setminus \langle B, t \rangle ^*$.  Then clearly $a \not \in A$, and, as in the proof of Proposition \ref{lower}, we can show that $A \cup \{a\}$ is a weakly $t$-independent set of size $m+1$ in $G$. 

{\bf Claim 2.}  Suppose that $A$ is a weakly $t$-independent set in $G$ of size $m$.   Then $$n
\geq \sum_{h=1}^{\lfloor t/2 \rfloor} \bc{m}{h}+1$$.

{\it Proof of Claim 2.}  Define $$\langle A, \lfloor t/2 \rfloor \rangle^* =\bigcup_{h=1}^{\lfloor t/2 \rfloor} h \star A.$$  By (\ref{onestar}), $0 \not \in \langle A, \lfloor t/2 \rfloor \rangle^* $, so we have $n-1 \geq |\langle A, \lfloor t/2 \rfloor \rangle^* |$.  Furthermore, by conditions (\ref{twostar}) and (\ref{threestar}), we see that $\langle A, \lfloor t/2 \rfloor \rangle^*$ has size exactly $$\sum_{h=1}^{\lfloor t/2 \rfloor} \bc{m}{h},$$ proving  our Claim.

To derive our upper and lower bounds for $w(G,t)$, we use the (rather crude) estimates that for positive integers $c$ and $d$ we have 
\begin{equation} \label{bin} \frac{(d+2-c)^c}{c!} \leq \bc{d+1}{c} \leq \sum_{h=0}^{c} \bc{d}{h} \leq \bc{d+c}{c} \leq \frac{(d+c)^c}{c!}.\end{equation}  Namely, using Claim 1 and (\ref{bin}) for $d=2m-2$ and $c=t$ we see that if $$n>\frac{(2m-2+t)^t}{t!},$$ then $G$ has a weakly $t$-independent set of size $m$, implying the lower bound $$w(G,t) \geq   \frac{(t!n)^{1/t}-t+2}{2}  -1 = \left( \frac{t!}{2^t} n \right)^{1/t}- \frac{t}{2}.$$  The upper bound for $w(G,t)$ follows similarly from Claim 2 and (\ref{bin}). $ \quad \Box$

Note that, for a fixed $t$, we have $$\liminf w(G,t) = \infty$$ as $|G|=n$ approaches $\infty$, in contrast to $$\liminf s(G,t) = 0$$ for each $t \geq 2$; in Section 1 we have seen that even $$\liminf \{s(G,t) | \mathrm{Ord}(G,t) \not = G\} = O(1)$$ as $t \geq 4.$   Thus $t$-independence and weak $t$-independence behave quite differently in abelian groups.

{\bf Acknowledgments.} The values of $s(\mathbb{Z}_n,4)$ and $s(\mathbb{Z}_n,5)$ for $n \leq 200$, on which Conjecture \ref{conj_t=4,5} was based, were computed by Nick Laza; we thank him for his time and efforts.

\end{document}